\documentclass[numreferences]{kluwer}
 \usepackage{amsfonts, amssymb, amsthm, euscript, amscd, latexsym,bm, mathrsfs}
\newtheorem{thm}{\sc Theorem}
\newtheorem{lem}{\sc Lemma}
\newtheorem{cor}{\sc Corollary}

 \newtheorem{Asn}{\sc Assumption}

\newcommand{\disp}{\displaystyle}

\newcommand{\eps}{\varepsilon}
\newcommand{\pl}{\partial}
\newcommand{\gt}{\geqslant}
\newcommand{\lt}{\leqslant}
\newcommand{\sub}{\subset}

\newcommand{\dl}{\delta}

 \newcommand{\Gm}{\Gamma}
 \newcommand{\Dl}{\Delta}
 \newcommand{\la}{\lambda}
 \newcommand{\lac}{\lambda_{\mathrm{c}}}
 \newcommand{\lastep}{\lambda_{\mathrm{step}}}
 \newcommand{\La}{\Lambda}
 \newcommand{\sg}{\sigma}
 \newcommand{\sgc}{\sigma_{\mathrm{c}}}
\newcommand{\dd}{\diagdown}
\newcommand{\om}{\omega}
\renewcommand{\mc}{\mathcal}
\newcommand{\Om}{\Omega}

\newcommand{\msr}{\mathscr}

\newcommand{\td}{\tilde}

\newcommand{\<}{\langle}
\renewcommand{\>}{\rangle}
\newcommand{\x}{\times}
\newcommand{\mto}{\mapsto}

\newcommand{\C}{{\rm C}}

\newcommand{\ind}{\mathbb I}

\newcommand{\skp}{\noalign{\vskip 0.5pc}}
\newcommand{\skps}{\noalign{\vskip 0.2pc}}
\newcommand{\text}{\textrm}

\newcommand{\fdot}{\,\cdot\,}
\def\Rnu{{\mathbb R}}

\def\Nnu{{\mathbb N}}

\def\ffi{\varphi}

\def\intl{\int\limits}

\def\com#1{}

\def\eee#1{ \begin{equation} #1 \end{equation} }
\def\aa#1{ \begin{eqnarray*} #1 \end{eqnarray*} }
\def\aaa#1{ \begin{eqnarray} #1 \end{eqnarray} }
\def\mm#1{ \begin{array}  #1 \end{array} }

\begin{document}
\begin{article}
\begin{opening}
\title{A mathematical approach to the nonequilibrium work theorem}
\author{Evelina \surname{Shamarova}\email{evelina@cii.fc.ul.pt}}
\institute{Grupo de F\'isica-Matem\'atica da Universidade de
Lisboa}

 \begin{abstract}
 We develop a mathematical approach to the
 nonequilibrium work theorem which is traditionally referred to in
 statistical mechanics as Jarzynski's identity.  We suggest
 a mathematically rigorous formulation and proof of the identity.
\end{abstract}

 \keywords{Nonequilibrium work theorem, Jarzynski's identity,
 nonequilibrium statistical mechanics,
 probability measures on phase space paths.}

 \classification{2000 Mathematics Classification}{60J25, 60J35, 82C05}

 \end{opening}

 \section{Introduction}
 The
  nonequilibrium work theorem is an equation in statistical mechanics that
  relates the free
  energy difference $\Dl F$ to the work $W$ carried out on a system
  during a nonequilibrium transformation.
  The identity appeared in different, but as we show below
  equivalent, formulations in  \cite{Jar97_1} and \cite{Jar97}
  in 1997, and in the series of papers \cite{bochkov1}--\cite{bochkov4} in 1977--1981.

  In physics
  literature, the identity is usually written in the form
  \eee{
  \label{jarzynski}
  \< e^{-\beta W} \> = e^{-\beta\Dl F},
  }
  where the average is taken over
  all possible system trajectories in the phase space,
  and $\beta$ is an inverse temperature.
  The identity first appeared in this form in
  \cite{Jar97_1} and \cite{Jar97}.
  Traditional equilibrium thermodynamics tells us that $\< W \> \gt \Dl F$ while
  the transformation of the system is infinitely slow. The identity
  (\ref{jarzynski})
  is a stronger statement, and in addition to this, it is valid
  for arbitrary transformations of the system.
  The identity is used effectively in computer simulations,
  as well as in experimental physics, to calculate the free energy difference
  between two states of the system by running many trajectories and taking
  the average value of $e^{-\beta W}$ (see \cite{Crooks2000}--\cite{Dellago2},
  \cite{Hummer_Szabo}-- \cite{liphardt} and references therein).

  The paper \cite{Jar07} discusses the connection between two different versions
  of the identity, and  shows that the papers \cite{bochkov1}--\cite{bochkov4}
  use a different definition of work.
  The identity obtained in \cite{bochkov1}--\cite{bochkov4} (referred to below as
   Bochkov--Kuzovlev's identity) reads:
   \[
   \< e^{-\beta W_0} \> = 1,
   \]
  where $W_0$ is the work (in Bochkov--Kuzovlev's sense) performed on the system,
  and the angle brackets have the same meaning as in (\ref{jarzynski}).
  The present paper shows that Jarzynski's and
  Bochkov--Kuzovlev's identities easily follow from each other.

 Since the identities involve taking an ``average over trajectories'', it
 is  natural to interpret this average as the expectation relative to a
 probability measure on trajectories, while assuming
 that the system evolves stochastically.
 In terms of expectations the identities can be represented by the formulas
 \[
 \mathbb E[e^{-\beta W}] = e^{-\beta\Dl F} \quad \text{and} \quad \mathbb E[e^{-\beta W_0}] = 1,
 \]
 where $\mathbb E$
 is the expectation relative to a probability measure on
 phase space paths.
 For this probability measure,
 some analytical assumptions
 under which the identities hold are found.

  \section{Notation and assumptions}
 Let us assume that the evolution of our system is described by
 a Markov process $\Gamma_t(\om)$, $t\in[0,T]$,  given through its transition
 density function.
  Let $X=\Rnu^{2d}$ be the phase space for our system,
  i.e.\ the set of values of $\Gamma_t$.
  We assume that at time $t=0$ the distribution
  on the phase space $X$ is given by the following density
  function:
  \[
  q_{\la_0}(x)=\frac{e^{-\beta \mc H(x,\la_0)}}{\int_X e^{-\beta \mc
  H(x',\la_0)}\, dx'}
  =\frac1{Z_{\la_0}}\: e^{-\beta \mc H(x,\la_0)},
  \]
  where $\mc H(\fdot, \la) : X \to \Rnu$ is a Hamiltonian
  parametrized by an externally controlled parameter $\la\in \La$,
  $\Lambda \sub \Rnu^l$ is an open
  set;
  $\beta = 1/(k_BT)$, $k_B$ is the Boltzmann constant,
  $T$ is the temperature of the system, $Z_{\la_0}=\int_X e^{-\beta \mc
  H(x,\la_0)}\, dx$ is the partition
  function.
  We assume that for all $\la_0 \in \Lambda$,
  $\int_X e^{-\beta \mc H(x,\la_0)}\, dx < \infty$.
  We consider the situation when the external parameter $\la$ is a function
  of time $[0,T]\to \Lambda$, i.e.\ we  actually consider a time-dependent
  Hamiltonian $\mc H(x,\la(t))$.
  Let $\mathbb E_{\la_0}$ denote the expectation
  relative to the measure $q_{\la_0}(x)dx$.
 Below we assume that the changing in time external parameter $\la$ belongs to the space
\[
V[0,T]=\{\la: [0,T]\to \Lambda, \, \la = \lac+\lastep; \; \lac \in \C^V, \,
\lastep\in \mc L_{\mathrm{step}}\}
\]
where $\C^V=\C^V[0,T]$ is the space
continuous function of bounded variation on $[0,T]$,
and
\[
\mc L_{\mathrm{step}} = \mc L_{\mathrm{step}}[0,T] = \{
 \la(\,\cdot\,)=
 \sum\nolimits_{i=0}^{n-1}\la_i\ind_{[t_{i},t_{i+1})}(\,\cdot\,) + \la_n
 \}
 \]
 is the space of right continuous step functions
 corresponding to different partitions
 $\mc P = \{0=t_0 < \cdots < t_n=T\}$ of $[0,T]$
 and different finite sets of values $\{\la_i\}$.

 To emphasize the fact that the function
 $\la\in V[0,T]$ gives rise to the process $\Gm_t$,
 we will use the notation $\Gm_t^\la$.
 Let $p_\la(s,x,t,y)$, $s,\,t\in [0,T]$, $s<t$, $x,\,y\in X$,
 be the transition density function for $\Gamma_t^\la$.
  Let $X^{[0,T]}$ denote the space of all paths $[0,T]\to X$.
  In terms of $p_\la$ we can construct a probability measure $\mathbb L_{\la}$
  on $X^{[0,T]}$ by means of the finite dimensional distributions
  \aaa{
  \disp
 &&\hspace{-3mm}\int_{X^{[0,T]}} f(\om(t_0), \om(t_1), \ldots,\om(t_n))\,\mathbb L_{\la}(d\om)
  =\int_X dx_0 \, q_{\la(t_0)}(x_0)
  \label{prob_dist}\\
  \skps
  \disp
  &&\hspace{-3mm}\intl_X  \hspace{-1mm} dx_1 \, p_\la(t_0, x_0, t_1, x_1)
  \cdots
   \intl_X \hspace{-1mm} dx_n \, p_\la(t_{n-1},x_{n-1}, t_n, x_n)\,f(x_0,x_1,\ldots,
   x_n), \nonumber
  }
  where
  $\{0 = t_0 < \cdots < t_n =T\}$ is a partition
  of the interval $[0,T]$, and $f : X^{n+1} \to \Rnu$ \ is
  a bounded and measurable function.
  By Kolmogorov's extension theorem,
  the right hand side of this equality defines
  a probability measure on the minimal $\sg$-algebra
  of $X^{[0,T]}$ generated by all cylindrical sets.
  We denote this $\sg$-algebra by $\sgc(X^{[0,T]})$.
  We assume that $p_\la(s,x,t,y)$ satisfies one of the assumptions
  (see \cite{skorokhod}, chapter 2, paragraph 1)
  that guaranties that the measure $\mathbb L_{\la}$ is concentrated on the right continuous
  trajectories without discontinuities of the second kind.
  Also, we assume that the $\sg$-algebra $\sgc(X^{[0,T]})$
  is augmented with all subsets of $\mathbb L_{\la}$-null sets.


  We take $X^{[0,T]}$ as the probability space,
  i.e. we set $\Om = X^{[0,T]}$.
  Then, $\mathbb L_{\la}$ is the distribution of $\Gm_{t}^\la$,
  and $\mathbb L_\la$-a.s., $\Gm^\la_t(\om) = \om(t)$.

 For each fixed $\bar \la \in \La$,
 we introduce another transition density function $p(s,x,t,y,\bar\la)$
 which represents the situation when the system evolves being controlled
 by a constant in time parameter.
 The following below Assumption \ref{asn1} is the key assumption under which
 the nonequilibrium work theorem holds.
 \begin{Asn}
 \label{asn1}
 If $\la|_{[s,t]}\equiv \bar\la$, then $p_\la(s,x,t,y)= p(s,x,t,y,\bar\la)$,
 where $p(s,x,t,y,\bar\la)$ conserves the canonical distribution
 on the phase space $X$. Specifically, it satisfies the identity
  \eee{
   \label{new1}
   \int_X q_{\bar\la}(x)\, p(s,x,t,y,\bar\la)\, dx   =   \, q_{\bar\la}(y).
   }
 \end{Asn}
 We make two further assumptions:
 \begin{Asn}
 \label{asn2}
 If $\la$ is constant on $[s,t)$,
 and discontinuous at the point $t$,
 then for Lebesgue almost all $x\in X$,
 \aa{
 &&\lim_{\dl\to 0+}
 \int_X  dy \, p(s,x, t -\dl, y, \la(s))
 \int_X dz \, p_\la(t- \dl, y, t, z) f(z) \\
 &&= \int_X  dy \, p(s,x, t, y, \la(s)) f(y),
  }
 where $f: X\to \Rnu$ is bounded and continuous.
 \end{Asn}
 \begin{Asn}
 \label{asn2'}
 For all $t\in [0,T)$, and for all compacts $K\sub X$,
 \[
 \lim_{\dl\to 0+} \sup_{x\in K}
 \Bigl[\int_X p(t,x,t+\dl,y,\la(t))f(y) dy - f(x)\Bigr] = 0,
 \]
 where $f: X\to \Rnu$ is bounded and continuous.
 \end{Asn}

   \begin{lem}
   \label{lem1}
   Let
   $\mc P =\{0=t_0<t_1< \cdots < t_n=T\}$
   be a partition of $[0,T]$, and let
   $
   \la(\,\cdot\,)=
   \sum_{i=0}^{n-1}\la_i\ind_{[t_{i},t_{i+1})}(\,\cdot\,) + \la_n
   $. We assume that if
   $\la|_{[s,t]}\equiv \bar \la$, where $\bar \la\in \La$ is constant,
    then
   $p_\la(s,x,t,y)= p(s,x,t,y,\bar\la)$.
   Further we assume
   that the function $p_\la$ satisfies Assumption \ref{asn2}.
   Then, for all $s,\,t$, $0\lt s < t \lt T$, for all $x\in X$,
   and for all bounded and continuous functions $f: X \to \Rnu$,
  \aaa{
   & &\hspace{-4mm} \int_X dy \, p_\la(s,x,t,y)\, f(y)= \int_X dx_1\,p(s,x,\tau_1, x_1,\la(s))\nonumber\\
    \label{9}
   & & \hspace{-4mm}\int_X dx_2\,  p(\tau_1,x_1,\tau_2,x_2,\la(\tau_1))
   \cdots
  \int_X dy \,   p(\tau_k,x_k,t,y,\la(\tau_k))\, f(y), \quad
  }
  where
  $\{s< \tau_1 < \cdots < \tau_k < t\} = (\mc P \cup \{s, t\}) \cap
  [s,t]$. Moreover, the right hand side of (\ref{9}) defines
  a transition density function.
   \end{lem}
  \begin{proof}
  First we prove that the right hand side of (\ref{9})
  defines a transition density function.
  We have to  verify the Chapman--Kolmogorov equation.

  Let $r\in (s,t)$, and $\td \mc P = \mc P \cup \{s,r,t\}$.
  Further let
 $\mc P_1 =\{s<\tau_1 < \tau_2 < \cdots < \tau_l\} = \td\mc P \cap [s,r]$,
 and
 $\mc P_2 = \{r < \xi_1 < \xi_2 < \cdots < \xi_m < t \}= \td\mc P \cap [r,t]$.
 Then, $\mc P_1 \cup \mc P_2  = \td \mc P$,
 and  $[\tau_l,\xi_1)$ is the
 interval of the partition $\mc P$ such that $r\in [\tau_l,\xi_1)$.
 We have:
 \eee{
 \label{12}
 \mm{
 {l}
 \disp
 \int_X dx_1\, p(s,x,\tau_1,x_1,\la(s))
 \int_X dx_2\, p(\tau_1,x_1,\tau_2, x_2,\la(\tau_1))
 \cdots \\
 \skp\disp
 \int_X dy \, p(\tau_l,x_l,r, y,\la(\tau_l)) \int_X dz_1\, p(r,y,\xi_1,z_1,\la(r)) \cdots\\
 \skp\disp
 \int_X dz_m \, p(\xi_{m-1},z_{m-1},\xi_m, z_m,\la(\xi_{m-1}))
  \, p(\xi_m,z_m,t,z,\la(\xi_m))\,.
 }
 }
 Note that $\la(r)=\la(\tau_l)$,
 and hence, by the Chapman--Kolmogorov equation for
 the transition density function $p(\fdot,\fdot,\fdot,\fdot,\la(\tau_l))$,
 we have:
 \[
 \int_X dy \, p(\tau_l,x_l,r,y,\la(\tau_l))\, p(r,y,\xi_1,z_1,\la(r))=
 p(\tau_l,x_l,\xi_1,z_1,\la(\tau_l)).
 \]
 Let us take a bounded and continuous function $f$, and show (\ref{9}).
 \eee{
 \label{long}
  \mm{
  {lll}
  \disp
  \int_X p_\la(s,x,t,y)f(y)dy  & = &
  \disp \int_X dx'_1 \, p(s,x,\tau_1-\dl_1, x'_1,\la(s))\\
  \skp
  & & \disp
  \int_X dx_1 \, p_\la(\tau_1-\dl_1,x'_1,\tau_1, x_1)
  \\
 \skp
 &\cdots& \disp \int_X  dx'_{k}\,
 p(\tau_{k},x_{k}, t-\dl_{k}, x'_{k}, \la(\tau_{k}))\\
 \skp
 & & \disp \int_X dy \, p_\la(t-\dl_{k},x'_{k}, t, y) f(y)\\
  }
  }
 which holds for all $\dl_1,\ldots, \dl_{k}$ smaller
 than the mesh of $\mc P$. Taking the repeated limit
 $\lim_{\dl_1\to 0} \lim_{\dl_2\to 0} \cdots \lim_{\dl_{k}\to 0}$
 of the right hand side of (\ref{long}), applying Lebesgue's theorem,
 and taking into consideration Assumption \ref{asn2}, we obtain
 the right hand side of (\ref{9}).
 \end{proof}
\section{Jarzynski's identity}
 Let $\la = \lac + \lastep$
 be the decomposition of $\la\in V[0,T]$ into a sum of a $\lac \in \C^V[0,T]$ and a $\lastep\in \mc L_{\mathrm{step}}[0,T]$.
 Further let $\pl_\la \mc H: X\x \La \to \Rnu^l$ denote
 the partial derivative with respect to the second argument
 (i.e.\ with respect to the control parameter).
 We assume that the restriction of
  $\pl_\la \mc H$ to $X\x \la([0,T])$ is bounded and measurable.
  Also, we assume that for each fixed $\bar \la \in \la([0,T])$,
  $\mc H(\fdot,\bar \la)$ is bounded and measurable.
 Everywhere below, the probability space $\Om$ is the space $X^{[0,T]}$.
 We define the work $W_\la: \Om \to \Rnu$
 performed on the system by the formula
 \eee{
 \label{8}
 \mm{
 {lll}
 \disp
 W_\la(\om) & = &
 \disp \int_0^T \bigl\langle \pl_\la \mc H(\Gamma_t^\la(\om),\la(t))\, ,\,d\lac(t)\bigr\rangle_{\Rnu^l} \\
  \skps\disp
 & + & \disp\sum_{i=1}^{n} \bigl(\mc H(\Gm^\la_{t_i}(\om), \la(t_i))
 - \mc H(\Gm^\la_{t_i}(\om),\la(t_i-0))\bigr),
 }
 }
  where  $\{0=t_0 < t_1 < \cdots < t_n = T\}$ are discontinuity points of $\la$,
  $\langle \; , \, \rangle_{\Rnu^l}$ is the scalar product in $\Rnu^l
  \supset\Lambda$,
  and the integral on the right hand side is the
  Lebesgue--Stieltjes integral, i.e.\ the sum
  of the Lebesgue--Stieltjes integrals with respect to the components of $\lac$.
  In the following, we skip the sign of the scalar product
  in the first term of~(\ref{8}), and simply write
 \aa{
 W_\la(\om) & = &
\int_0^T \pl_\la \mc H(\Gamma_t^\la(\om),\la(t))\, d\lac(t) \\
 & + & \sum_{i=1}^{n} \bigl(\mc H(\Gm^\la_{t_i}(\om), \la(t_i))
 - \mc H(\Gm^\la_{t_i}(\om),\la(t_i-0))\bigr).
 }
  Let $F_\la = -\frac1{\beta}\, \ln
  Z_\la$ (free energy of the system),
  and let $\Dl F = F_{\la(T)} - F_{\la(0)}$
  (free energy difference).
  Let $\mathbb E_{\mathbb L_{\la}}$ denote the expectation
  relative to the measure $\mathbb L_\la$.
 \subsection{Jarzynski's Identity, case $\la\in \mc L_{\mathrm{step}}[0,T]$}
  Clearly, if
  $\la = \sum_{i=0}^{n-1} \la_i \ind_{[t_{i},t_{i+1})} + \la_n \in \mc L_{\mathrm{step}}[0,T]$,
  \[
  W_\la(\om)=
   \sum_{i=1}^{n} \bigl(\mc H(\Gm^\la_{t_i}(\om), \la_{i})
 - \mc H(\Gm^\la_{t_i}(\om),\la_{i-1})\bigr).
 \]
 \begin{thm}[Jarzynski's identity: case ${\la\in\mc L_{\mathrm{step}}[0,T]}$]
 \label{Jarzynski_theorem}
 Let $\la\in \mc L_{\mathrm{step}}[0,T]$, and let
 the transition density function $p_\la$ satisfy
 Assumptions \ref{asn1} and \ref{asn2}.
 Further let
 $\mc H(\fdot, \bar\la)$ be bounded and measurable on $X$ for each fixed $\bar \la \in \la([0,T])$.
 Then the function $e^{-\beta W_\la}$ is $\mathbb L_\la$-integrable,
 and
 \[
  \mathbb E_{\mathbb L_\la}[e^{-\beta W_\la}] = e^{-\beta\Dl F}.
 \]
 \end{thm}

  \begin{proof}[Proof]
   $\mathbb L_\la$-a.s.,
   \[
   W_\la(\om) =
   \sum_{i=1}^{n} \Bigl(\mc H(\om(t_i),\la_{i}) - \mc
   H(\om(t_i),\la_{i-1})\Bigr).
   \]
   Without loss of generality, we can assume that
   $t_0=0$, $t_n=T$, adding these points with zero jumps
   if necessary.
 Note that $W_\la(\om)$ is a cylinder function.
 By Lemma \ref{lem1},
 \eee{
 \label{long000}
  \mm{
  {lll}
  \disp
  &&\mathbb E_{\mathbb L_\la}[e^{-\beta W_\la}]
   = \disp \int_{X^{[0,T]}} \, e^{-\beta \sum_{i=1}^{n}
  (\mc H(\om(t_i),\la_{i}) - \mc
  H(\om(t_i),\la_{i-1}))}\, \mathbb L_\la(d\om)\\
  \skp
  & = &  \disp
  \int_X dx_0
  \, q_{\la_0}(x_0)
  \int_X dx_1 \, p(t_0,x_0,t_1, x_1,\la_0)\\
  \skp
 & \cdots &\disp \int_X  dx_{n-1}\, p(t_{n-2},x_{n-2},t_{n-1},
 x_{n-1},\la_{n-2})\\
 \skp
 & & \disp \int_X dx_n \, p(t_{n-1},x_{n-1}, t_{n},x_{n},\la_{n-1})
 \,
  e^{-\beta \sum_{i=1}^{n}
  (\mc H(x_i,\la_{i}) - \mc H(x_i,\la_{i-1}))}.
  }
  }
  Note that
  Assumption \ref{asn1} implies that for all $s<t$,
  for all $y\in X$, and $\bar\la\in \Lambda$,
  \eee{
  \label{chaing_simplify_argument}
  \int_X   dx \: \frac{e^{-\beta \mc H(x,\bar\la)}}
        {e^{-\beta \mc H(y,\bar\la)}}
         \; p(s,x,t,y,\bar\la) = 1.
  }
  Taking into account this, and
  changing the order of integration in~(\ref{long000}),
  we obtain
  \eee{
   \label{long1}
   \mm{
   {l}
   \mathbb E_{\mathbb L_\la}[e^{-\beta W_\la}]
    = \disp
    \frac1{Z_{\la_0}}\,
   \int_X dx_n \,
   e^{-\beta \mc H(x_{n}, \la_n)} \\
   \skp
   \disp
   \int_X dx_{n-1}
   \frac
    {e^{-\beta \mc H(x_{n-1}, \la_{n-1})}}
    {e^{-\beta \mc H(x_n,\la_{n-1})}} \:
    p(t_{n-1},x_{n-1},t_{n},x_{n},\la_{n-1})\\
    \skp
    \cdots\\
    \skp
    \disp \int_X  dx_0 \,
    \frac
    {e^{-\beta \mc H(x_0, \la_0)}}
    {e^{-\beta \mc H(x_1,\la_0)}} \;
    p(t_0,x_{0},t_1,x_{1},\la_{0}).
   }
   }
  Starting from the end, we replace
  each integral in  (\ref{long1}) with $1$,
  which is valid by the relation (\ref{chaing_simplify_argument}),
  until we reach the very first integral
  (taken with respect to $x_n$), which we replace with $Z_{\la_n}$.
  Noticing that $\la_0 = \la(0) $ and $\la_n = \la(T)$,
  we obtain
  \[
  \mathbb E_{\mathbb L_\la}[e^{-\beta W_\la}] =
  \frac{Z_{\la(T)}}{Z_{\la(0)}}.
  \]
  The theorem is proved.
 \end{proof}
 \subsection{Jarzynski's Identity, case: {$\la \in \C^V[0,T]$}}
 \begin{thm}[Jarzynski's identity, case: {$\la \in \C^V[0,T]$}]
 \label{Jarzynski_theorem_gen}
 Let $\la\in \C^V[0,T]$,
 and let the transition density function $p_\la$
 satisfy  Assumptions \ref{asn1}, \ref{asn2}, and \ref{asn2'},
 and the probability distribution $\mathbb L_\la$ of
 $\Gm^\la_t$ be given by  (\ref{prob_dist}).
 In addition, let the following assumptions
 be fulfilled:
 \begin{enumerate}
 \setcounter{enumi}{3}
 \item
 \label{asn3}
 If $\la^n\in \mc L_{\mathrm{step}}[0,T]$, and as $n\to\infty$,
  $\la^n \to \la$ uniformly on $[0,T]$, then
 $\mathbb L_{\la^n} \to \mathbb L_{\la}$ weakly relative
 to the family of bounded continuous cylinder functions;
 \item
\label{asn4'}
 The function $\mc H(\fdot, \bar\la)$ is bounded and measurable on
 $X$ for each fixed $\bar \la \in \la([0,T])$;
 \item
 \label{asn4}
 The function $\pl_\la \mc H$ is bounded on $X\x \la \bigl([0,T]\bigr)$;
 \item
 \label{asn44}
 The functions $\pl_\la \mc H(x,\fdot)$
 are equicontinuous as a family of functions parameterized by $x\in X$;
 \item
 \label{asn5}
 The function $\pl_\la \mc H(\fdot,\bar \la)$ is continuous for each fixed
 $\bar\la \in \Lambda$.

 \end{enumerate}
 Then, the function
 $e^{-\beta W_\la}$ is $\mathbb L_\la$-integrable,
 and
 \eee{
  \label{J_identity}
  \mathbb E_{\mathbb L_\la}[e^{-\beta W_\la}] = e^{-\beta\Dl F}.
 }
 \end{thm}
\begin{lem}
 Let $\la\in \C^V[0,T]$. Then,
  under Assumptions \ref{asn4}--\ref{asn5},
  the function $W_\la(\om)$ is $\sg_c$-measurable.
  \end{lem}
  \begin{proof}
  Below, for an arbitrary small $\eps$ we construct a
  $\sg_c$-measurable function $F: \Om \to \Rnu$ such that
  \eee{
  \label{verify}
  \sup_\Om |F(\om) - W_\la(\om)| < \eps.
  }
  We find a function $\la_{\mathrm{step}}=
  \sum_{i=0}^{n-1} \la(t_i) \ind_{[t_i,t_{i+1})}+ \la(t_n) \in \mc L_{\mathrm{step}}[0,T]$
  such that for all $\om\in \Om$, for all $t\in [0,T]$,
  \[
  |\pl_\la \mc H(\Gamma_t^\la(\om),\la(t))- \pl_\la \mc H(\Gamma_t^\la(\om),\la_{\mathrm{step}}(t))|<\eps.
  \]
  This is possible by Assumption \ref{asn44}.
  Thus, it suffices to find a measurable function $F$
  verifying (\ref{verify}) for
  $W_\la$ of the form
  $\int_s^r \pl_\la \mc H(\Gamma_t^\la(\om),\bar \la)\, d\lac(t)$,
  where $\bar \la \in \la([0,T])$ is fixed, and $\la_c$ is not a constant on $[s,r]$.

  By assumption, $\mathbb L_\la$-a.s., the paths of $\Gm^\la_t$ are
  right continuous and without discontinuities of the second kind.
  This implies that the map $\Gm^\la: [0,T] \x \Om \to X$
  is $(\msr B([0,T])\otimes \sg_c, \msr B(X))$-measurable
  (where $\msr B([0,T])$ and $\msr B(X)$ are Borel $\sg$-algebras)
  which follows from \cite{skorokhod} (Chapter 2, Theorem 11).
  By Assumption \ref{asn4}, $\pl_\la \mc H(\fdot,\fdot)$ is bounded,
  say by a constant $M$. We divide the ball of radius $M$ in $\Rnu^l$
  with the center in the origin,
  by a finite number of sets $O_i$ whose diameter is smaller than $\eps/V_0^T[\la_c]$,
  where $V_0^T[\la_c]$ denotes the variation of $\la_c$ on $[0,T]$.
  Further let $A_i = \pl_\la \mc H(\fdot,\bar \la)^{-1}(O_i)$,
  and $C_i = (\Gm^\la)^{-1}(A_i)\sub [0,T]\x \Om$.
  The sets $A_i$ are open by Assumption \ref{asn5}.
  The sets $C_i$ are $\msr B([0,T])\otimes \sg_c$-measurable
  by the $(\msr B([0,T])\otimes \sg_c, \msr B(X))$-measurability
  of the map $\Gm^\la$. Fix $x_i\in A_i$, and consider the function
  $\Phi(\om,t) = \sum_i \pl_\la \mc H(x_i,\bar \la)\ind_{C_i}(\om,t)$.
  Clearly,
  $
  \sup_{t\in [0,T],\om\in \Om} |\pl_\la \mc H(\om(t),\bar \la) - \Phi(\om,t)|<\eps / V_0^T[\la_c].
  $
  We define
  \[
  F(\om) = \int_s^r \Phi(\om,t) \la(dt)= \sum_i \pl_\la \mc H(x_i,\bar \la)
  \mu_\la(\{t: (\om,t) \in C_i\}),
  \]
   where $\mu_\la$ is the Lebesgue-Stieltjes measure on $[0,T]$
   corresponding to the function $\la$. Fubini's theorem
   implies that the function $\om\mto \mu_\la(\{t: (\om,t) \in C_i\})$
   is $\sg_c$-measurable.  The inequality (\ref{verify}) is obviously satisfied.

  \end{proof}

\begin{proof}[Proof of Theorem \ref{Jarzynski_theorem_gen}]
We take a sequence of partitions
 $\mc P_n = \{0=t^n_0 < t^n_1 < \cdots < t^n_n = T\}$, and
  consider the functions
 \[
 \la^n(t) =
 \sum_{i=0}^{n-1}\la(t^n_i)\ind_{[t^n_{i},t^n_{i+1})}(t)
 + \la(T).
 \]
 Clearly, $\la_n \rightrightarrows \la$ on $[0,T]$,
 and $V_0^T[\la^n] < V_0^T[\la]$,
 where $V_0^T$ is the variation on $[0,T]$.
 We set $\la^n_i=\la(t^n_i)$, and define the functions
 \[
 \ffi_n(\om) = e^{-\beta \sum_{i=1}^{n}
 (\mc H(\om(t^n_i),\la^n_{i}) - \mc H(\om(t^n_i),\la^n_{i-1}))}.
 \]
 By Theorem \ref{Jarzynski_theorem},
 \[
 \int_{\Om} \ffi_n(\om)\, \mathbb L_{\la^n} (d\om) =
 \frac{Z_{\la(T)}}{Z_{\la(0)}}
 \]
 for all $n$. We denote $\ffi(\om) = e^{-\beta W_\la(\om)}$, and prove that
  \aaa{
  \label{term1}
  & &\lim_{n\to\infty}\int_\Om
   \bigr(\ffi_n(\om) -  \ffi(\om) \bigl)\mathbb L_{\la^n}(d\om)
   = 0, \\
   \label{term2}
  & &\lim_{n\to\infty}\int_\Om
   \bigr(\ffi_n(\om) -  \ffi(\om) \bigl)\mathbb L_\la(d\om)
   = 0,\\
   \label{term3}
  & &\lim_{n\to\infty}\lim_{m\to\infty}\int_\Om
   \bigr(\ffi_n(\om) -  \ffi(\om) \bigl)\mathbb L_{\la^m}(d\om)
   = 0.
  }
  For this, we first replace the functions $\ffi_n$
  in (\ref{term1})--(\ref{term3}) with
  more suitable functions $\hat \ffi_n$ such that
  $\bigl(\ffi_n(\om) - \hat\ffi_n(\om)\bigr)\to 0$
  as $n\to\infty$, uniformly on $\Om$. We have:
 \eee{
 \label{2t}
 \mm{
 {ll}
 \disp
 & \disp \hspace{-4mm}\sum_{i=1}^{n}
 (\mc H(\om(t^n_i),\la^n_{i}) -
 \mc H(\om(t^n_i),\la^n_{i-1}))
 - \int_0^T \pl_\la \mc H(\om(t),\la^n(t)) \, d\la^n(t)\nonumber\\
 &  \disp\hspace{-4mm} = \sum_{i=1}^{n}
  \pl_\la \mc H(\om(t^n_i),\td\la^n_i)
  (\la^n_{i} - \la^n_{i-1})
 - \sum_{i=1}^{n}
  \pl_\la \mc H(\om(t^n_i), \la^n_i)
  (\la^n_{i} - \la^n_{i-1}),
 }
 }
 where in the first term in (\ref{2t}) we applied the mean value
 theorem to each summand, and chose $\td\la_i^n \in [\la_i^n,\la_{i+1}^n]$.
 Since by assumption $\mc \pl_\la H(x,\fdot)$ is
 equicontinuous, the absolute value of the difference
 in (\ref{2t}) does not exceed $\eps V_0^T[\la]$ where
 $\eps$ is chosen so that $|\pl_\la\mc H(x,\bar\la_1) - \pl_\la\mc H(x,\bar\la_2)|<\eps$
 whenever $|\bar \la_1 - \bar \la_2| < \dl$, and $\dl$ is chosen
 by the equicontinuity argument.
 The relation (\ref{2t}) shows that if we
 prove  (\ref{term1})--(\ref{term3})
 with $\td \ffi_n(\om) = e^{-\beta \int_0^T \pl_\la \mc H(\om(t),\la^n(t)) \, d\la^n(t)}$
 substituted for $\ffi_n$, then we prove
 (\ref{term1})--(\ref{term3}).
 We define the functions:
 \aaa{
 \label{om_n}
 &\om^n(t,\om) =  \sum\nolimits_{i=0}^{n-1} \om(t^n_i) \ind_{[t^n_{i},t^n_{i+1})}(t)
 +\om(T);
 \nonumber\\
 \label{td-om_n}
 &\td \om^n(t,\om) = \om(0) + \sum\nolimits_{i=1}^{n} \om(t^n_i) \ind_{(t^n_{i-1},t^n_{i}]}(t); \\
 & \td \la^n(t)  = \la(0) +
 \sum_{i=1}^{n}\la(t^n_i)\ind_{(t^n_{i-1},t^n_{i}]}(t). \nonumber
 }
 With the help of these functions, the second term in (\ref{2t})
 can be represented as
  \eee{
   \label{est1}
  \mm{
  {lll}
  \disp
  \int_0^T \pl_\la \mc  H(\om(t),\la^n(t))d\la^n(t)
  & = & \disp\int_0^T \pl_\la \mc  H(\om^n(t,\om),\la^n(t))d\la^n(t) \\
  & = & \disp\int_0^T \pl_\la \mc  H(\td \om^n(t,\om),\td \la^n(t))d\la(t).
  }
  }
  By  Assumption \ref{asn44},
  we choose an $\eps > 0$ so that
    $|\pl_\la\mc H(\td \om^n(t,\om),\td \la^n(t)) - \pl_\la\mc H(\td \om^n(t,\om),\la(t))|<\eps$
  whenever $\sup_{t\in [0,T]}|\td \la^n(t) - \la(t)| < \dl$, and $\dl$ is chosen
  by the equicontinuity argument.
  This means that the relations (\ref{term1}), (\ref{term2}), and  (\ref{term3})
  are equivalent to
  \eee{
  \label{term1aaa}
   \mm{
   {l}
  \disp\lim_{n\to\infty}\int_\Om
   \bigr(\hat \ffi_n(\om) -  \ffi(\om) \bigl)\mathbb L_{\la^n}(d\om)
   = 0, \\
   \skp
  \disp\lim_{n\to\infty}\int_\Om
   \bigr(\hat\ffi_n(\om) -  \ffi(\om) \bigl)\mathbb L_{\la}(d\om) = 0, \\
   \skp
   \disp\lim_{n\to\infty} \lim_{m\to\infty}\int_\Om
   \bigr(\hat\ffi_n(\om) -  \ffi(\om) \bigl)\mathbb L_{\la^m}(d\om) = 0,
  }
  }
  where
  $\hat \ffi_n(\om) = e^{-\beta \int_0^T \pl_\la \mc  H(\td \om^n(t,\om),\la(t))d\la(t)}$,
  and $\td \om^n(t,\om)$ is given by (\ref{td-om_n}).
  We show that the relations (\ref{term1aaa}) hold.
  Since $\pl_\la \mc H(\fdot,\fdot)$ is bounded by Assumption \ref{asn4},
  $\la$ is a function of bounded variation on $[0,T]$,  and the exponent is Lipschitz
  on bounded domains, for all $m$ we obtain the estimate
  \aaa{
  && \hspace{-10mm} \left|
  \int_\Om
  \bigr(\hat\ffi_n(\om) -  \ffi(\om) \bigl)\mathbb L_{\la^m}(d\om)
  \right|\label{show}\\
  && \hspace{-10mm} < K_L
  \int_\Om \mathbb L_{\la^m}(d\om) \int_0^T |\la|(dt)
  \bigr|\pl_\la \mc H(\td \om^n(t,\om),\la(t))- \pl_\la \mc H(\om(t),\la(t))
  \bigl| \nonumber\\
  && \hspace{-10mm}= K_L \int_0^T |\la|(dt)
  \int_\Om \mathbb L_{\la^m} (d\om)
  \,|\pl_\la\mc H(\td \om^n(t,\om),\la(t))-\pl_\la\mc H(\om(t),\la(t))|, \nonumber
  }
 where $K_L$ is the Lipschitz constant for the exponent, $|\la|(dt)$
 is the Lebesgue-Stieltjes measure corresponding to the total
 variation function $|\la|(t)$,
 and  Fubini's theorem has been applied to pass to the last integral.
 The same estimate holds for $\mathbb L_\la$. Namely,
 \aaa{
 && \hspace{-10mm} \left|
  \int_\Om
  \bigr(\hat\ffi_n(\om) -  \ffi(\om) \bigl)\mathbb L_{\la}(d\om)
  \right|\label{show1}\\
  && \hspace{-10mm}= K_L \int_0^T |\la|(dt)
  \int_\Om \mathbb L_{\la} (d\om)
  \,|\pl_\la\mc H(\td \om^n(t,\om),\la(t))-\pl_\la\mc H(\om(t),\la(t))|.
 \nonumber
 }
 We would like to show that
 \aaa{
 \label{repeated}
 & &\hspace{-2mm}\lim_{n\to\infty}\!
 \disp\int_\Om \mathbb L_{\la^n}(d\om)
   \bigr|\pl_\la \mc H(\td \om^n(t,\om),\la(t))- \pl_\la \mc H(\om(t),\la(t))
   \bigl| = 0, \qquad \\
 \label{single}
 & &\lim_{n\to\infty}\!
 \disp\int_\Om \mathbb L_\la(d\om)
   \bigr|\pl_\la \mc H(\td\om^n(t,\om),\la(t))- \pl_\la \mc H(\om(t),\la(t))
   \bigl| = 0, \qquad\\
 \label{repeated1}
 \lim_{n\to\infty}&&\hspace{-3mm}\lim_{m\to\infty}\!
 \disp\int_\Om \mathbb L_{\la^m}(d\om)
   \bigr|\pl_\la \mc H(\td\om^n(t,\om),\la(t))- \pl_\la \mc H(\om(t),\la(t))
   \bigl| = 0. \qquad
 }
 Let
 $t\in (t^n_{i-1},t^n_{i}]$,
 and let $\dl_n = t_i^n -t$.
 Then $\td \om^n(t,\om) = \om (t + \dl_n)$.
 Let us show (\ref{repeated}).
 We have:
 \aa{
  && \int_{\Om}
  \mathbb L_{\la^n} (d\om)
  \,|\pl_\la\mc H(\om(t+\dl_n),\la(t))-\pl_\la\mc H(\om(t),\la(t))|\nonumber\\
    = & & \int_X q_{\la_0}(x_0) dx_0 \int_X  dx \, p_{\la^n}(t_0,x_0,t,x) \nonumber\\
  & & \int_X dx_i\, p(t,x,t+\dl_n,x_i,\la(t)) \, |\pl_\la\mc H(x,\la(t))-\pl_\la\mc H(x_i,\la(t))|.
  \quad
  }
  Note that
  \aa{
  F_{n} (x)= \int_X  dx_i \, p(t,x,t + \dl_n,x_i, \la(t)) \, \Bigl|\pl_\la\mc H(x_i,\la(t))
  -\pl_\la\mc H(x,\la(t))\Bigr|
  }
  converges to zero uniformly in $x$ running over compacts in $X$
  which follows from Assumption \ref{asn2'}.
  Indeed,
  \aa{
  |F_n(x)|^2 \lt \int_X \hspace{-1mm} dx_i\, p(t,x,t + \dl_n,x_i, \la(t)) \, \Bigl|\pl_\la\mc H(x_i,\la(t))
  -\pl_\la\mc H(x,\la(t))\Bigr|^2.
  }
  The right hand side of this inequality converges to zero
  uniformly in $x\in K \sub X$, where $K$ is an arbitrary compact.
  This easily follows
  from Assumption \ref{asn2'} after we separate in the right hand side
  the terms depending on $x$ and on $x_i$.
  Define the measures:
  \aa{
  &\mu_{\la^n}(A)  =  \int_X q_{\la_0}(x_0)\, dx_0 \int_A p_{\la^n}(t_0,x_0,t,x) \, dx,\\
  &\mu_\la(A) = \int_X q_{\la_0}(x_0) \, dx_0 \int_A p_{\la}(t_0,x_0,t,x)\, dx.
  }
  By Assumption \ref{asn3}, $\mathbb L_{\la^n}$ converges weakly
  to $\mathbb L_{\la}$ relative to the family of
  bounded continuous cylinder function,
  whenever $\la^n \rightrightarrows \la$.
  This implies that $\mu_{\la^n} \to \mu_\la$ weakly (relative
  to the family of bounded continuous functions), as $\la^n \rightrightarrows \la$.
  By Prokhorov's theorem, the family  $\{\mu_{\la^n}, \mu_\la\}$
  of probability measures on $X$ is tight. We fix an arbitrary $\eps > 0$
  and find a compact $K_\eps$ such that $\mu_{\la^n}(X \dd K_\eps) < \eps$
  for all $n$.
  Since $F_n \rightrightarrows 0$ on $K_\eps$, and all
  $F_n$
  are bounded on $X$ by a constant, say $M$, not depending on $n$,
  we can find
  an $N \in \Nnu$ such that $|F_n(x)|<\eps$
  on $K_\eps$ for all $n>N$. We obtain the estimate:
  \[
   \Bigl|
   \int_X \mu_{\la^n}(dx) F_n(x)
   \Bigr|
   \lt  \sup_{x\in K_\eps} |F_n(x)| + M \mu_{\la^n}(X\dd K_\eps)
   \lt   \eps + M \eps.
  \]
  This proves (\ref{repeated}).
  The relation (\ref{single}) follows from the
  right continuity of $\om$, Assumption \ref{asn5}, and Lebesgue's theorem.
  The relation (\ref{repeated1}) follows from
  Assumption \ref{asn3}, and from (\ref{single}).
  Application of Lebesgue's theorem to the integral
  (\ref{show}) (taken over $[0,T]$ with respect to $|\la|(dt)$)
  implies
  that as $n\to\infty $, the integral (\ref{show}) converges to zero.
  This implies (\ref{term1aaa}), and thus (\ref{term1}), (\ref{term2}), and
  (\ref{term3}) are proved.
  By (\ref{term1}), we obtain:
  \aa{
   &0& = \lim_{n\to \infty} \int (\ffi_n \mathbb L_{\la^n} - \ffi \mathbb L_{\la^n})
   = \lim_{n\to \infty} \int \ffi_n \mathbb L_{\la^n} - \lim_{n\to\infty}
  \int \ffi \mathbb L_{\la^n}\\
  & & = \frac{Z_{\la(T)}}{Z_{\la(0)}} - \lim_{n\to\infty}
  \int \ffi \mathbb L_{\la^n}.
  }
  On the other hand, (\ref{term2}) and (\ref{term3}) imply:
  \aa{
  &0& = \lim_{n\to \infty}\lim_{m\to\infty} \int (\ffi_{\la^n}
  \mathbb L_{\la^m} - \ffi \mathbb L_{\la^m}) \\
  & & = \lim_{n\to\infty}\lim_{m\to\infty} \int \ffi_n \mathbb L_{\la^m}
  - \lim_{m\to\infty} \int \ffi \mathbb L_{\la^m} \\
  & & =  \lim_{n\to\infty} \int \ffi_n \mathbb L_\la - \lim_{m\to\infty} \int \ffi \mathbb L_{\la^m}
  = \int \ffi \mathbb L_\la - \lim_{m\to\infty} \int \ffi \mathbb L_{\la^m}.
  }
  Comparing the last two relations gives:
  \[
  \int_\Om \ffi(\om) \mathbb L_\la(d\om) = \frac{Z_{\la(T)}}{Z_{\la(0)}}.
  \]
  The theorem is proved for the case $\la\in \C^V[0,T]$.
  \end{proof}

  \subsection{Jarzynski's identity for $\la \in V[0,T]$ and its corollaries}
 \begin{cor}[Corollary of Theorem \ref{Jarzynski_theorem_gen}]
 \label{corollary1}
 Let the assumptions of Theorem~\ref{Jarzynski_theorem_gen}
 be fulfilled, and let
 $f:~X\to \Rnu$ be bounded and continuous.
 Then,
  \eee{
  \label{corollary_relation}
  \mathbb E_{\mathbb L_\la}[(f\circ\pi_T) \, e^{-\beta \, W_\la}]
  = \mathbb E_{\la(T)}[f]\; \mathbb E_{\mathbb L_\la}[e^{-\beta \, W_\la}],
  }
  where $\pi_t: X^{[0,T]}\to X, \om \mto \om(t)$ is the
  evaluation mapping,
  $\mathbb E_{\la(T)}$ is the expectation relative to the measure
  $\frac1{Z_{\la(T)}}\, e^{-\beta \,\mc H(x,\la(T))}dx$.
 \end{cor}

 \begin{proof}
  Assuming that $\la\in \mc L_{\mathrm{step}}[0,T]$, we repeat the
  argument
  of~(\ref{long}) and~(\ref{long1}), while
  using the
  relation~(\ref{chaing_simplify_argument}).
  Specifically, we obtain:
  \aa{
  & &\mathbb E_{\mathbb L_\la}[(f\circ\pi_T) \, e^{-\beta \,
  W_\la}]=
    \frac1 {Z_{\la(t_0)}}
    \intl_X
    {e^{-\beta \mc H(x_{n}, \la(t_{n}))}}
    f(x_n) \, dx_n\\
   & & \,
   \intl_X
    \frac
    {e^{-\beta \mc H(x_{n-1}, \la(t_{n-1}))}}
    {e^{-\beta \mc H(x_{n},\la(t_{n-1}))}} \:
    p(t_{n-1},x_{n-1},t_{n},x_{n},\la(t_{n-1})) \, dx_{n-1}
    \ldots
    \\
  \skp
  & & \cdots\\
  \skp
  & & \intl_X
    \frac
    {e^{-\beta \mc H(x_0, \la(t_0))}}
    {e^{-\beta \mc H(x_1,\la(t_0))}} \;
    p(t_0,x_{0},t_1,x_{1},\la(t_{0})) \, dx_{0} \\
   \skps
   & =&
   \mathbb E_{\la(T)}[f] \; \frac{Z_{\la(t_n)}}{Z_{\la(t_0)}}
   = \mathbb E_{\la(T)}[f] \; \mathbb E_{\mathbb L_\la}[e^{-\beta \,
    W_\la}]\,.
   }
   Now let $\la \in \C^V[0,T]$, and let the constant $M_f$ be
   such that $\sup_{x\in X}|f(x)| < M_f$.
   As in the proof of Theorem \ref{Jarzynski_theorem_gen},
   we find a sequence of step functions $\la^n$
   converging to $\la$ uniformly on $[0,T]$.
   The identity (\ref{corollary_relation}) holds for each $\la^n$.
   The argument of passing to the limit as $n\to \infty$
   is similar to the argument used
   in the proof of Theorem \ref{Jarzynski_theorem_gen}.
   We define the function $F: \Om \to \Rnu$, $F(\om) = f(\om(T))$, and
   repeat the argument of Theorem \ref{Jarzynski_theorem_gen}
   until the inequality (\ref{show1})
   with the following replacements:
   $\ffi \leftrightarrow F \ffi$, $\ffi_n \leftrightarrow F \ffi_n$,
   $\td \ffi_n \leftrightarrow F \td\ffi_n$, $\hat\ffi_n \leftrightarrow F \hat\ffi_n$.
   The right hand sides of (\ref{show}) and (\ref{show1}) will
   be the same as in the proof of Theorem \ref{Jarzynski_theorem_gen}
   but the Lipschitz constant $K_L$ will be replaced with $K_L M_f$.
   The part of the proof of Theorem \ref{Jarzynski_theorem_gen}
   following after the inequality (\ref{show1}) remains unchanged
   until the last two arguments of passing to the limit.
   Those arguments now will be:
   \aa{
   &0& = \lim_{n\to \infty} \int (F\ffi_n \mathbb L_{\la^n} - F\ffi \mathbb L_{\la^n})
   = \lim_{n\to \infty} \int F \ffi_n \mathbb L_{\la^n} - \lim_{n\to\infty}
  \int F \ffi \mathbb L_{\la^n}\\
  & & =  \mathbb E_{\la(T)}[f]\; \mathbb E_{\mathbb L_\la}[e^{-\beta \, W_\la}]
   - \lim_{n\to\infty}\int F \ffi \mathbb L_{\la^n}.
  }
  On the other hand,
  \aaa{
  &0& = \lim_{n\to \infty}\lim_{m\to\infty} \int (F\ffi_{\la^n}
  \mathbb L_{\la^m} - F\ffi \mathbb L_{\la^m}) \nonumber\\
  & & = \lim_{n\to\infty}\lim_{m\to\infty} \int F\ffi_n \mathbb L_{\la^m}
  - \lim_{m\to\infty} \int F\ffi \mathbb L_{\la^m} \nonumber \\
  & & = \int F\ffi \mathbb L_\la - \lim_{m\to\infty} \int F\ffi \mathbb L_{\la^m}.
  \quad \label{98}
  }
  This implies the identity (\ref{corollary_relation}).
  \end{proof}
  \begin{cor}[Corollary of Theorem \ref{Jarzynski_theorem_gen}]
  \label{corollary22}
  Let the assumptions of Corollary \ref{corollary1}
  except Assumption \ref{asn3} of Theorem \ref{Jarzynski_theorem_gen}
  be fulfilled with respect to the function $\la\in \C^V[0,T]$.
  Let $a\in \Lambda$, and let $\la^a = \la \,\ind_{[0,T)} + a \,\ind_{\{T\}}$.
  Let us assume that
  \vspace{-2mm}
  \begin{itemize}
  \item[]
  If $\la^n\in \mc L_{\mathrm{step}}[0,T]$, and
  $\la^n \rightrightarrows \la^a$  on $[0,T]$, then
  $\mathbb L_{\la^n} \to \mathbb L_{\la^a}$ weakly relative
  to the family of bounded continuous cylinder functions;
  \end{itemize}
  Then,
  \eee{
  \label{cor_id_modified}
  \mathbb E_{\mathbb L_{\la^a}}[(f\circ\pi_T) \, e^{-\beta \, W_\la}]
  = \mathbb E_{\la(T)}[f]\; \frac{Z_{\la(T)}}{Z_{\la(0)}}.
  }
  \end{cor}
  \begin{proof}
  If $\la^a \in \mc L_{\mathrm{step}}[0,T]$, i.e.
  $\la^a = \sum_{i=1}^{n-1} \la_i \,\ind_{[t_i,t_{i+1})} + a\,\ind_{\{T\}}$,
  then the representation
  (\ref{9}) shows that $\mathbb E_{\mathbb L_{\la^a}}[(f\circ\pi_T) \, e^{-\beta \, W_{\la}}]$
  is the same for the functions $\la^a$ with different $a$.
  Hence,
  \[
  \mathbb E_{\mathbb L_{\la^a}}[(f\circ\pi_T) \, e^{-\beta \, W_{\la}}]
  = \mathbb E_{\la(T)}[f]\; \frac{Z_{\la(T)}}{Z_{\la(0)}}\,.
  \]
  Note that in the relations (\ref{show}) and (\ref{show1}) we integrate with respect
  to $|\la|(dt)$ where $|\la|$ is continuous.
  Hence, by Lebesgue's theorem,
  the investigation of convergence at the point $T$
  in the integrals with respect
  to $\mathbb L_{\la_m}(d\om)$ and $\mathbb L_{\la}(d\om)$
  is not necessary.
  Finally we note that by assumption, and by the relations (\ref{show1}) and
  (\ref{single}), with $\la^a$ substituted for $\la$, for the first term  in the identity (\ref{98})
  we obtain:
  \[
  \lim_{n\to\infty}\lim_{m\to\infty} \int F\ffi_n \mathbb L_{\la^m}
  = \int F\ffi \, \mathbb L_{\la^a}.
  \]
  This proves (\ref{cor_id_modified}).
  \end{proof}

 Let $[s,t]\sub [0,T]$ be an arbitrary subinterval,
 $X^{[s,t]}$ be the space of all functions
 on $[s,t]$. We introduce further notation:
 Let
 $\C^V[s,t]$ be the space of continuous functions of bounded variation on $[s,t]$,
 $\mc L_{\mathrm{step}}[s,t]$ be the space of right continuous step functions
 on $[s,t]$ taking a finite number of values.
 Below, on $X^{[s,t]}$ we define the distributions
 $\mathbb L^x_{\la;s,t}$ and $\mathbb L_{\la;s,t}$,
 $\la \in V[0,T]$, $x\in X$.
 Let $s < \tau_1 < \cdots < \tau_{k-1} < t$ be a partition of $[s,t]$,
 $f: X^{k+1}\to \Rnu$ be bounded and measurable. The finite dimensional
 distributions of $\mathbb L_{\la;s,t}^x$ are defined by the formula:
  \[
  \mm{
  {l}
  \disp
 \int_{X^{[s,t]}} f(\om(s), \om(\tau_1), \ldots, \om(t))\,\mathbb L_{\la;s,t}^x(d\om)
  =
  \\
  \skps
  \disp
  \intl_X  \hspace{-1mm} dx_1 \, p_\la(s, x, \tau_1, x_1))
  \cdots
   \intl_X \hspace{-1mm} dx_k \, p_\la(\tau_{k-1},x_{k-1}, t, x_k)\,f(x, x_1,\ldots, x_k).
  }
  \]
  We extend $\mathbb L^x_{\la;s,t}$ to $\sg_c(X^{[s,t]})$ by Kolmogorov's extension
  theorem, and notice that under assumptions that we made,
  this measure it concentrated
  on the right continuous paths without discontinuities of the second kind.
   Let $f$ be $\mathbb L_{\la;s,t}^x$-integrable. We define
  \aaa{
  \label{expectations_formula}
   \mathbb L_{\la;s,t}[f] & =
  \mathbb E_{\la(s)}\bigl[\mathbb L^{\bullet}_{\la;s,t}[f]\bigr].
  }

 \begin{thm}[Jarzynski's identity, case: {$\la \in V[0,T]$}]
 \label{Jarzynski_theorem_gen_gen}
 Let $\la\in V[0,T]$, and let $\{t_1 < \dots < t_k = T\}$ be
 the set of its discontinuity points. Further let
 the transition density function $p_\la$
 of the Markov process $\Gm^\la_t$ satisfy
 Assumptions \ref{asn1}, \ref{asn2}, and \ref{asn2'},
 and the probability distribution $\mathbb L_\la$ be given by (\ref{prob_dist}).
 Let Assumptions \ref{asn4'}, \ref{asn4}, \ref{asn44}, and \ref{asn5}
 of Theorem \ref{Jarzynski_theorem_gen} be fulfilled.
 Additionally, we assume that
 \begin{enumerate}
 \setcounter{enumi}{7}
 \item
 \label{asn8}
 If $\la^m \in \mc L_{\mathrm{step}}[0,T]$, and
 $\la^m \rightrightarrows \la$ on $[0,T]$, as $m\to\infty$,
 then for each $i$, $0 \lt i < n$,
 $\mathbb L_{\la^m; t_i,t_{i+1}} \to \mathbb L_{\la; t_i,t_{i+1}}$ weakly relative
 to the family of bounded continuous cylinder functions on $X^{[t_i,t_{i+1}]}$;
 \item
 \label{asn9}
 For all bounded and continuous functions $f: X^{[t_i,t_{i+1}]} \to \Rnu$,
 the function $X \to \Rnu, \; x\mto \mathbb L^x_{\la;t_i,t_{i+1}}[f]$
 is continuous.
 \end{enumerate}
 Then, the function
 $e^{-\beta W_\la}$ is $\mathbb L_\la$-integrable,
 and
 \eee{
  \label{J_identity_gen}
  \mathbb E_{\mathbb L_\la}[e^{-\beta W_\la}] = e^{-\beta\Dl F}.
 }
 \end{thm}
 The continuity of functions $X^{[t_i,t_{i+1}]} \to \Rnu$ is understood with respect to the topology
 of pointwise convergence.
 \begin{proof}[Proof of Theorem \ref{Jarzynski_theorem_gen_gen}, case {$\la \in V[0,T]$}]
 Let $\la = \lac + \lastep$, where $\lac \in \C^V[0,T]$, and
  $\lastep = \sum_{i=0}^{n-1} \la_i \ind_{[t_{i},t_{i+1})}+ \la_n \in \mc L_{\mathrm{step}}[0,T]$.
  On each interval $[t_{i-1},t_i]$ we define the function
  $\om_i = \om|_{[t_{i-1},t_i]}$,
  and identify each $\om \in X^{[0,T]}$ with the sequence $(\om_1, \ldots, \om_n)$.
  Since $\Gm^\la_t$ is a Markov process, we obtain:
  \[
  \mm{
  {l}
  \disp
  \int_{X^{[0,T]}} e^{-\beta W_\la(\om)} \mathbb L_{\la}(d\om)
  =
  \int_X dx_0 \, q_{\la(t_0)}(x_0)
  \int_{X^{[t_0,t_1]}} \mathbb L_{\la;t_0,t_1}^{x_0}(d\om_1)\\
  \skps
  \disp
  \int_{X^{[t_1,t_2]} }\mathbb L_{\la;t_1,t_2}^{\om_1(t_1)}(d\om_2)
  \cdots
  \int_{X^{[t_{n-1},t_n]}} \mathbb L_{\la;t_{n-1},t_n}^{\om_{n-1}(t_{n-1})}(d\om_n)
  e^{-\beta  W_\la(\om_1, \ldots, \om_n)}
  }
  \]
 where
 \[
 \mm{
 {ll}
 \disp
 W_\la(\om_1, \ldots, \om_n) & = \disp  \sum_{i=1}^{n} \int_{t_{i-1}}^{t_{i}}
 \pl_\la \mc H(\om_i(t),\la(t))d\lac(t)\\
 \skps
 & \disp + \sum_{i=1}^{n}
 \bigl(\mc H(\om_i(t_i),\la(t_i)) - \mc H(\om_i(t_i),\la(t_i-0))\bigr).
 }
 \]
 On each interval $[t_{i-1},t_i]$, we
 define continuous functions $\la_{(i)}(t) = \la(t)$, $t\in [t_{i-1},t_{i})$,
 and $\la_{(i)}(t_{i})= \la(t_{i}-0)$, $1\lt i \lt n$, and introduce a notation:
 \[
 W_{\la_{(i)}}(\om_i)  =  \int_{t_{i-1}}^{t_i} \pl_\la \mc H(\om_i(t),\la(t))d\la_{(i)}(t).
 \]
 We obtain:
 \[
 \mm{
 {l}
 \disp
 \intl_{X^{[0,T]}} e^{-\beta W_\la(\om)} \mathbb L_{\la}(d\om)
  = \disp \int_X dx_0 \, q_{\la(t_0)}(x_0)\\
 \skps
  \disp  \intl_{X^{[t_0,t_1]}}
  e^{-\beta W_{\la_{(1)}}(\om_1)}
  e^{-\beta \bigl(
  \mc H(\om_1(t_{1}),\la(t_{1}))   - \mc H(\om_1(t_{1}), \la(t_{1}-0))
  \bigr)}
   \mathbb L_{\la;t_0,t_1}^{x_0}(d\om_1)\\
  \skps
  \disp  \intl_{X^{[t_1,t_2]}}
  e^{-\beta W_{\la_{(2)}}(\om_2)}
  e^{-\beta \bigl(
  \mc H(\om_2(t_{2}),\la(t_{2}))   - \mc H(\om_2(t_{2}), \la(t_{2}-0))
  \bigr)}
   \mathbb L_{\la;t_1,t_2}^{\om_1(t_1)}(d\om_2) \cdots \\
 \skps
  \disp
  \intl_{X^{[t_{n-1},t_{n}]}} \hspace{-5mm}
  e^{-\beta W_{\la_{(n)}}(\om_{n})}
  e^{-\beta \bigl(
  \mc H(\om_n(t_{n}),\la(t_{n}))   - \mc H(\om_n(t_{n}), \la(t_{n}-0))
  \bigr)}
  \mathbb L_{\la;t_{n-1},t_n}^{\om_{n-1}(t_{n-1})}(d\om_{n}).
  }
 \]
 Let $\mathbb E_{\bar \la}$ denote the expectation relative to the measure
 $q_{\bar\la}(x)dx$.
 We can rewrite the above relation in terms of expectations:
 \aaa{
  \label{above}
 &&\hspace{-5mm}\mathbb E_{\mathbb L_\la} [e^{-\beta W_\la}]\\
 && \hspace{-5mm} = \mathbb E_{\la(t_0)}[
  \mathbb E_{\mathbb L^{\bullet}_{\la;t_0,t_1}}\,[e^{-\beta W_{\la_{(1)}}(\om_1)}
 e^{-\beta(\mc H(\om_1(t_{1}),\la(t_{1}))   - \mc H(\om_1(t_1), \la(t_1-0)))}
 \mathbb E_{\mathbb L^{\om_1(t_1)}_{\la;t_1,t_2}} \nonumber\\
&& \hspace{-5mm}[\cdots \mathbb E_{\mathbb L^{\om_{n-1}(t_{n-1})}_{\la;t_{n-1},t_n}}
 \bigl[e^{-\beta W_{\la_{(n)}}(\om_n)} e^{-\beta \bigl(
  \mc H(\om_n(t_{n}),\la(t_{n}))   - \mc H(\om_n(t_{n}), \la(t_{n}-0))
  \bigr)} \bigr]
 \cdots \, ]]]. \nonumber
 }
 Let $F: X^{[t_i,t_{i+1}]} \to \Rnu$  be bounded and continuous.
 We obtain:
 \aaa{
 && \mathbb E_{\la(t_i-0)} \Bigl[
  e^{-\beta ( \mc H(\bullet, \la(t_i)) -\mc H(\bullet, \la(t_i-0))}\,
 \mathbb E_{\mathbb L^{\bullet}_{\la;t_i,t_{i+1}}}[F] \,
 \Bigr]\nonumber\\
 \label{55}
  & & = \frac{Z_{\la(t_i)}}{Z_{\la(t_i-0)}}\,
  \mathbb E_{\la(t_i)}\Bigl[
  \mathbb E_{\mathbb L^{\bullet}_{\la; t_i,t_{i+1}}}[F]
  \Bigr]
  = \frac{Z_{\la(t_i)}}{Z_{\la(t_i-0)}}\,
 \mathbb E_{\mathbb L_{\la;t_i,t_{i+1}}}[F].
 \quad
 }
 Note that in (\ref{above}) $F$ is always
 a function of the form
 \eee{
 \label{measurable}
 X^{[t_i,t_{i+1}]}\to \Rnu, \;
 \om \mto  e^{-\beta W_{\la_{(i+1)}}(\om)}f(\om(t_{i+1})),
 }
 where $f: X\to \Rnu$ is bounded and continuous by Assumptions \ref{asn5}
 and \ref{asn9}.
 Taking
 the expectation $\mathbb E_{\mathbb L_{\la; t_i,t_{i+1}}}$
 of the function (\ref{measurable}), and applying Corollary \ref{corollary22}
 we obtain:
 \eee{
 \label{this_id}
 \mathbb E_{\mathbb L_{\la; t_i,t_{i+1}}}[ e^{-\beta W_{\la_{(i+1)}}(\om)}f(\om(t_{i+1}))]
  = \frac{Z_{\la(t_{i+1}-0)}}{Z_{\la(t_i)}}\, \mathbb E_{\la(t_{i+1}-0)}[f].
 }
 Applying  the identities (\ref{55}) and (\ref{this_id}) to (\ref{above}), we obtain:
 \eee{
 \label{last}
 \mathbb E_{\mathbb L_\la} [e^{-\beta W_\la}]
 = \frac{Z_{\la(t_1-0)}}{Z_{\la(t_0)}}\, \frac{Z_{\la(t_1)}}{Z_{\la(t_1-0)}}
 \cdots \frac{Z_{\la(t_{n}-0)}}{Z_{\la(t_{n-1})}}\,\frac{Z_{\la(t_n)}}{Z_{\la(t_{n}-0)}}
 = \frac{Z_{\la(T)}}{Z_{\la(0)}}.
 }
 \end{proof}
\begin{cor}[Corollary of Theorem \ref{Jarzynski_theorem_gen_gen}]
 \label{corollary2}
 Let the assumptions of Theorem~\ref{Jarzynski_theorem_gen_gen}
 be fulfilled, and let
 $f:~X\to \Rnu$ be bounded and continuous.
 Then,
  \eee{
  \label{corollary_relation_gen}
  \mathbb E_{\mathbb L_\la}[(f\circ\pi_T) \, e^{-\beta \, W_\la}]
  = \mathbb E_{\la(T)}[f]\; \mathbb E_{\mathbb L_\la}[e^{-\beta \, W_\la}].
  }
\end{cor}
\begin{proof}
 We repeat the argument that we used in the proof of Theorem \ref{Jarzynski_theorem_gen_gen}
 to obtain (\ref{above}) in connection to the expression
 $\mathbb E_{\mathbb L_\la}[(f\circ\pi_T) \, e^{-\beta \, W_\la}]$.
 Instead of the very last expectation in (\ref{above}), we obtain
 \[
 \mathbb E_{\mathbb L^{\om_{n-1}(t_{n-1})}_{\la; t_{n-1},t_n}}[ e^{-\beta W_{\la_{(n)}}(\om_n)}\, f(\om_n(t_n))
  \, e^{-\beta \bigl( \mc H(\om_n(t_n),\la(t_n))   - \mc H(\om_n(t_n), \la(t_n-0))\bigr)}].
 \]
 Applying Corollary \ref{corollary22}, we obtain:
 \aaa{
  \label{29}
 & & \mathbb E_{\mathbb L_{\la; t_{n-1},t_n}}[ e^{-\beta W_{\la_{(n)}}(\om)}f(\om(t_n))\,
 e^{-\beta \bigl( \mc H(\om_n(t_n),\la(t_n))   - \mc H(\om_n(t_n), \la(t_n-0))\bigr)}] \nonumber\\
 & & = \frac{Z_{\la(t_n-0)}}{Z_{\la(t_{n-1})}}\,
 \mathbb E_{\la(t_n-0)} [e^{-\beta(\mc H(\fdot,\la(t_n))   - \mc H(\fdot, \la(t_n-0)))}\, f] \nonumber\\
 & & =  \frac{Z_{\la(t_n-0)}}{Z_{\la(t_{n-1})}}\,
  \frac{Z_{\la(t_n)}}{Z_{\la(t_n-0)}}\, \mathbb E_{\la(t_n)}[f].
 }
 Instead of the multiplier
 $\frac{Z_{\la(t_{n}-0)}}{Z_{\la(t_{n-1})}}\,\frac{Z_{\la(t_n)}}{Z_{\la(t_{n}-0)}}$
 in the relation (\ref{last})
 we obtain $\frac{Z_{\la(t_n-0)}}{Z_{\la(t_{n-1})}}\,
  \frac{Z_{\la(t_n)}}{Z_{\la(t_n-0)}}\, \mathbb E_{\la(t_n)}[f]$
 according to (\ref{29}). This gives (\ref{corollary_relation_gen}).
\end{proof}

 \section{Bochkov--Kuzovlev's identity}
  Here we give a rigorous mathematical proof of the identity announced
  in \cite{bochkov1}--\cite{bochkov4} where the authors used a different definition
  of work to that in the papers \cite{Jar97_1}--\cite{Jar97}. The difference
  between the two definitions of work was analyzed in the paper
  \cite{Jar07}, and it was found that
   \[
    W(\om)-W_0(\om) = \mc H(\om(T),\la(T)) - \mc H(\om(T),\la(0)),
   \]
  where $W_0$ is the work in the Bochkov--Kuzovlev sense. We use this equality
  as the definition of the new work $W_0$, and will prove the following theorem:
  \begin{thm}
  Let the assumptions of Theorem \ref{Jarzynski_theorem_gen_gen}
  be fulfilled. Then,
  \eee{
  \label{BK_identity}
  \mathbb E_{\mathbb L_\la} [ e^{-\beta W_0}] = 1.
  }
  Moreover, the identities (\ref{BK_identity}) and (\ref{J_identity_gen})
  are equivalent.
  \end{thm}
  \begin{proof}
  Applying Corollary \ref{corollary2} and the identity
  (\ref{J_identity_gen})
  we obtain:
  \aa{
  \mathbb E_{\mathbb L_\la} [ e^{-\beta W_0}]
  & = & \mathbb E_{\mathbb L_\la}[ e^{-\beta W}]\,
  \mathbb E_{\la(T)} [e^{-\beta\bigl(\mc H(\om(T),\la(T)) - \mc H(\om(T),\la(0))\bigr)}]\\
  & = & \mathbb E_{\mathbb L_\la}[ e^{-\beta W}]\, \frac{Z_{\la(0)}}{Z_{\la(T)}}
  = \frac{Z_{\la(T)}}{Z_{\la(0)}}\,  \frac{Z_{\la(0)}}{Z_{\la(T)}} = 1.
   }
   This relation also shows that the identities
   (\ref{BK_identity}) and (\ref{J_identity_gen})
  are equivalent.
  \end{proof}

  \begin{acknowledgements}
  The author thanks Professor C.\ Dellago for attracting her
  attention to Jarzynski's identity, and Professor O.\ G.\ Smolyanov
  for useful discussions.
  This work was supported by the research grant
  of the Erwin Schr\"odinger Institute for mathematical physics,
  by the Austrian Science Fund (FWF) under START-prize-grant Y328,
  and by the research grant of the
  Portuguese Foundation for Science and Technology (FCT).
  \end{acknowledgements}


\end{article}

\end{document}